\documentclass[12pt]{amsart}
\usepackage {amsmath, amscd}
\usepackage {amssymb}
\usepackage {epsfig}
\usepackage {bm}
\usepackage {indentfirst} %indent the first par after section

\usepackage[active]{srcltx}

\advance\hoffset by -1.5cm

\numberwithin{equation}{section}

\def\<{\langle}
\def\>{\rangle}

\def\L{{\bf L}}

\def\BB{{\mathcal B}}

\def\HH{{\mathcal H}}

\def\JJ{{\mathcal J}}

\def\LL{{\mathcal L}}

\def\bbR{\mathbb{R}}

\newtheorem{lemma}{Lemma}[section]
\newtheorem{proposition}[lemma]{Proposition}%[section]
\newtheorem{theorem}[lemma]{Theorem}%[section]
\newtheorem{corollary}[lemma]{Corollary}

\theoremstyle{definition}
\newtheorem{remark}[lemma]{Remark}

\title{Mind Duggal Transforms}
\author{ C. Benhida}
\address{Chafiq Benhida\\
Universit\'e Lille 1\\
 Laboratoire Paul Painlev\'e \\
  UFR de Math\' ematiques, UMR 8524, { B\^at M2},
 59655 Villeneuve d'Ascq\\ France}\email{chafiq.benhida@math.univ-lille1.fr}

\begin{document}

%\maketitle

\begin{abstract}
It is known that if an operator $T$ is complex symmetric then its Aluthge transform is also complex symmetric. This Note is devoted to showing that  the Duggal transform doesn't inherit this property. For instance, we'll show that the Duggal transform isn't always complex symmetric
when $T$ is, as it was claimed in \cite{Ga}. 
\end{abstract}

\maketitle \footnotetext{2010 Mathematics Subject Classification;
Primary 47A05, 47A10, 47A11.
\par  Keywords; Aluthge  transform; complex symmetric operator; Duggal transform.
}

\maketitle

\section { Introduction }

Let ${\LL}(\HH)$ be the algebra of all bounded linear operators on a separable complex Hilbert space ${\HH}$.  For an operator
$T\in{\mathcal{L}}(\mathcal{H})$, $T^{\ast}$ denotes the adjoint of $T$. An
operator $T\in{\mathcal{L}}(\mathcal{H})$ is said to be \textit{normal} if
$T^{\ast}T=TT^{\ast}$, \textit{quasinormal} if $T^{\ast}T$ and $T$ commute,
\textit{binormal} if $T^{\ast}T$ and $TT^{\ast}$ commute, \textit{subnormal}
if there exists a Hilbert space $\mathcal{K}$ containing $\mathcal{H}$ and a
normal operator $N$ on $\mathcal{K}$ such that $N\mathcal{H}\subset
\mathcal{H}$ and $T=N|_{\mathcal{H}}$, and \textit{hyponormal} if $T^{\ast
}T-TT^{\ast}\geq0$. 

A {\it conjugation} on $\HH$ is an antilinear operator $C: {\HH}\rightarrow {\HH}$ which satisfies $\langle Cx, Cy \rangle=\langle y, x\rangle$ for
all $x,y\in {\HH}$ and $C^{2}=I$. An operator $T\in{\LL(\HH)}$ is said to be {\it complex symmetric} if
there exists a conjugation $C$ on ${\HH}$ such that $T= CT^{\ast}C$. 
%and \it{skew complex symmetric} \em if there exists a conjugation $C$ on
%$\HH$ such that $CTC = -T^*$. 
Many standard operators such as normal operators, algebraic operators of order $2$, Hankel matrices, finite Toeplitz
matrices, all truncated Toeplitz operators, and  Volterra integration operators are included in the class of complex symmetric operators.
Several authors have studied the structure of complex symmetric operators (see \cite{Ga 3}-\cite{GW}, \cite{JKLL}, and \cite{JKL} for more details). For spectral properties, see  also
\cite{BCKL}.

\medskip

Recall that for a given operator $T\in{\mathcal{L}}(\mathcal{H})$, we have the following writing 
$T=U|T|$ called the polar decomposition of $T$ where $U$ is a partial isometry with ($\ker U= \ker T$) and $|T|:= (T^*T)^{\frac{1}{2}}$.
 The Aluthge tranform of $T$ is the operator  $\widetilde{T}:=| T|^{\frac{1}{2}} U | T|^{\frac{1}{2}}$. This transform is playing an important role 
 in many aspects around the study of $T$ (see for example \cite{A}, \cite{BZ}, \cite{BZ1} and \cite{JKP}). An other operator connected to $T$ is the Duggal transform $T^D:=| T| U$ and will be considered in this paper
 concerning particularly complex symmetricity.

\section{ What happens for Duggal transform?}

We start by recalling the following result \cite[Theorem 3.1]{LZ2}:

\begin{proposition}\label{Zhu} If $T= \sum_{i=1}^{n-1} \lambda_i e_{i}\otimes e_{i+1}$  and $\lambda_i\ne 0$ for all $i$, then
 $T$ is complex symmetric if and only if $ |\lambda_i|= |\lambda_{n-i}|$ for evey $1\leq i\leq n-1$.
\end{proposition}

We'll show the following result which is an immediate consequence.

\begin{corollary}\label{Duggal} If $T= \sum_{i=1}^{n-1} \lambda_i e_{i}\otimes e_{i+1}$  and $\lambda_i\ne 0$ for all $i$, then
 its Duggal transform $T^D$ is complex symmetric if and only if $ |\lambda_i|= |\lambda_{n-1-i}|$ for evey $1\leq i\leq n-2$ .
\end{corollary}

{\it Proof.} One may without loss of generality assume that $\lambda_i >0$ for every $i $. Or equivalently, consider
  $T= \sum_{i=1}^{n-1} \lambda_i e_{i}\otimes e_{i+1}$  and $\lambda_i\ne 0$ for all $i$, 
 then 
$T= \sum_{i=1}^{n-1} |\lambda_i | f_{i}\otimes f_{i+1}$ where $f_1=e_1$ and $f_{i+1}:= \frac{  \bar{\lambda}_1\dots\bar{\lambda}_i }{|\lambda_1 
|\dots|\lambda_i |}e_{i+1}$ and of course
$\BB= \{f_{i},\  1\leq i\leq n \}$ is an orthonormal basis. If we write (and we'll do so for all matrices in the sequel) the matrix of $T$  
according to the basis $\BB$, we have

\[
T\cong{ \rm Mat} (T,\BB)=\begin{pmatrix}
0 & |\lambda_{1}| & 0 & & \ldots&  0\\
0 & 0 & |\lambda_{2}| & 0 & \ldots & 0\\
\vdots &\vdots &\ddots &\ddots & \ldots& 0\\
\vdots & \vdots & . & 0 & \ddots &  0\\
.& . & . & . & 0 & |\lambda_{n-1}|\\
0 & 0 & . & . &\ldots & 0 \\
\end{pmatrix}
\]

It has been shown \cite[Theorem 3.1]{LZ2} that $T$ is complex symmetric if and only if $|\lambda_{i}|=|\lambda_{n-i}|$ for every $1\leq i\leq n-1$.

A simple calculations shows that

\[
U= \begin{pmatrix}
0 & 1 & 0 & & \ldots&  0\\
0 & 0 & 1 & 0 & \ldots & 0\\
\vdots &\vdots &\ddots &\ddots & \ldots& 0\\
\vdots & \vdots & . & 0 & \ddots &  0\\
.& . & . & . & 0 & 1\\
0 & 0 & . & . &\ldots & 0 \\
\end{pmatrix}
\]

\[
|T|= \begin{pmatrix}
0 &   &  &  &  \\
 & |\lambda_{1}| &  &  &  \\
 &  &  &  \ddots &  \\
 &  &  &  &  |\lambda_{n-1}| \\
\end{pmatrix}
\]

Thus the Duggal transform is given by

\[
T^D= |T|U=\begin{pmatrix}
0 & 0 & 0 & & \ldots&  0\\
0 & 0 & |\lambda_{1} |& 0 & \ldots & 0\\
\vdots &\vdots &\ddots &\ddots & \ldots& 0\\
\vdots & \vdots & . & 0 & \ddots &  0\\
.& . & . & . & 0 & |\lambda_{n-2}|\\
0 & 0 & . & . &\ldots & 0 \\
\end{pmatrix}=0\oplus
\begin{pmatrix}
0 & |\lambda_{1} |& 0 & \ldots & 0\\
\vdots &\ddots &\ddots & \ldots& 0\\
 \vdots & . & 0 & \ddots &  0\\
 . & . & . & 0 & |\lambda_{n-2}|\\
 0 & . & . &\ldots & 0 \\
\end{pmatrix}
\]

and has (more or less) the same shape as $T$.

Using Proposition \ref{Zhu} and \cite[Lemma 1]{GW1} (which says that $A$ is complex symmetric if and only if $0\oplus A$ is complex symmetric), $T^D$ is complex symmetric if and only if $|\lambda_{i}|=|\lambda_{n-1-i}|$ for every $1\leq i\leq n-2$.

\begin{corollary}\label{Duggal1} If $T= \sum_{i=1}^{n-1} \lambda_i e_{i}\otimes e_{i+1}$  and $\lambda_i\ne 0$ for all $i$, then $T$ and
 its Duggal transforms $T^D$ are both complex symmetric if and only if $ |\lambda_1|= |\lambda_2|=\dots= |\lambda_{n-1}|$.
\end{corollary}

From what has been shown above, one easily infer that $T$  and $T^D$ are both complex symmetric if and only if $|\lambda_{1}|=|\lambda_{2}|=\ldots=|\lambda_{n-1}|$ 

which means that 

\[
T= \alpha\begin{pmatrix}
0 & 1 & 0 & & \ldots&  0\\
0 & 0 & 1 & 0 & \ldots & 0\\
\vdots &\vdots &\ddots &\ddots & \ldots& 0\\
\vdots & \vdots & . & 0 & \ddots &  0\\
.& . & . & . & 0 & 1\\
0 & 0 & . & . &\ldots & 0 \\
\end{pmatrix}
,\]

where $\alpha$ is arbitrary in $\bbR$.

Notice also that in this case all generalized Aluthge transforms of $T$ are complex symmetric 
with  the conjugaison $C(z_1,z_2,\ldots, z_n)= (\bar z_n, \ldots, \bar z_2,\bar z_1)$.

Also, we aren't facing the trivial case of a fixed point of Aluthge transform map which means that $T$ is not quasinormal and even more (see below).

Indeed, if 

\[
T= \begin{pmatrix}
0 & 1 & 0 & & \ldots&  0\\
0 & 0 & 1 & 0 & \ldots & 0\\
\vdots &\vdots &\ddots &\ddots & \ldots& 0\\
\vdots & \vdots & . & 0 & \ddots &  0\\
.& . & . & . & 0 & 1\\
0 & 0 & . & . &\ldots & 0 \\
\end{pmatrix}
,\]

then

\[
|T|= \begin{pmatrix}
0 &   &  &  &  \\
 & 1 &  &  &  \\
 &  &  &  \ddots &  \\
 &  &  &  &  1\\
\end{pmatrix}
\]

and

\[
T|T|= \begin{pmatrix}
0 & 1 & 0 & & \ldots&  0\\
0 & 0 & 1 & 0 & \ldots & 0\\
\vdots &\vdots &\ddots &\ddots & \ldots& 0\\
\vdots & \vdots & . & 0 & \ddots &  0\\
.& . & . & . & 0 & 1\\
0 & 0 & . & . &\ldots & 0 \\
\end{pmatrix}=T
\]

while

\[
|T|T= \begin{pmatrix}
0 & 0 & 0 &  \ldots&  0\\
0 & 0 & 1 &  \ldots & 0\\
\vdots &\vdots &\ddots  & \ddots& 0\\
\vdots & \vdots & \vdots & 0 & 1\\
0 & 0 &\ldots& 0 & 0 \\
\end{pmatrix}.
\]

Thus $T|T|\neq |T|T$ and $T$ is not quasinormal.

On the other hand, remark that this operator is binormal.

\section{Binormal operators and the symmetric property}

Recall that in ${\\L}(\HH)$, two operators $A$ and $B$  commute if $[A,B]=: AB-BA=0$

An operator $T$ in ${\\L}(\HH)$ is quasinormal if $T$ commutes with $T^*T$
and is said to be binormal if $TT^*$ commutes with $T^*T$ (or equivalently $[|T|, |T^*|]=0$).

\medskip

\begin{theorem}\label{bin} If $T= \sum_{i=1}^{n-1} \lambda_i e_{i}\otimes e_{i+1}$  and $\lambda_i\ne 0$ for all $i$, then
 $T$ is a binormal operator.
\end{theorem}

{\it Proof.} 

It's easy  to see that

\[
|T|= \begin{pmatrix}
0 &   &  &  &  \\
 & |\lambda_{1}| &  &  &  \\
 &  & |\lambda_{2}|  &  &  \\
 &  &  &  \ddots &  \\
 &  &  &  &  |\lambda_{n-1}| \\
\end{pmatrix}
\ \ \ {\text and } \ \ 
|T^*|= \begin{pmatrix}
|\lambda_{1}| &   &  &  &  \\
 & |\lambda_{2}| &  &  &  \\
 &  &  &  \ddots &  \\
 &  &  &  &  |\lambda_{n-1}| \\
 &  &  &  &  & 0\\
\end{pmatrix}
\]

\medskip

\[
|T||T^*|= \begin{pmatrix}
0 &   &  &  &  \\
 & |\lambda_{1}\lambda_{2}  | &  &  &  \\
 &  &    \ddots & & \\
 &  &  &  &  |\lambda_{n-2} \lambda_{n-1}| \\
  &  &  &  &  & 0\\
\end{pmatrix} = |T||T^*|
\]

\begin{remark}\label{R1}

\begin{enumerate}

 \item The claim in \cite[Proposition 3.1]{HMT} saying that: "a binormal operator  $T$ is
complex symmetric if and only if its Duggal transform is complex symmetric"
 is not true. One may construct easy examples from what has been shown above.

Indeed, it's enough to take

\[
T= \begin{pmatrix}
0 &  1 & 0 &  0 \\
0 &  0 & 2 &  0 \\
0 &  0 & 0 &  1 \\
0 &  0& 0 &  0 
 \end{pmatrix}
\]

 \item In the same paper, the authors are using as a fact that if an operator $T$ is  binormal and
complex symmetric with the polar decomposition $T=U|T|$ then $U$ is unitary. Also this claim is not true as one may see
from all our examples.
 
\end{enumerate}
\end{remark}

\section{Generalized Aluthge transforms and the symmetric property}

 Let $T\in{\LL(\HH)}$ with the polar decomposition $T=U |T|$. The Generalized Aluthge transform of $T$ is 
 the operator 
$\widetilde {T}(t)=|T|^{t}U |T|^{1-t}$ for $t\in [0,1]$.

\smallskip

One may see the following result as a generalization of the one given in section 2 (see also
\cite[Section 3]{LZ2}).

\begin{theorem}\label{gAlu}
If $T= \sum_{i=1}^{n-1} \lambda_i e_{i}\otimes e_{i+1}$  and $\lambda_i\ne 0$ for all $i$, then
 its  generalized Aluthge transform  $\widetilde {T}(t)$, for $t\in ]0,1]$.
 is complex symmetric if and only if $ |\lambda_i|^{t}|\lambda_{i+1}|^{1-t}= |\lambda_{n-1-i}|^{t}  |\lambda_{n-i}|^{1-t}$ for evey $1\leq i\leq n-2$ .
 In particular
 
 \begin{enumerate}
 \item 
Its Aluthge transform  $\widetilde {T}=\widetilde {T}(\frac{1}{2}) $ is complex symmetric if and only if $ |\lambda_i\lambda_{i+1}|= |\lambda_{n-1-i}\lambda_{n-i}|$ 
for evey $1\leq i\leq n-2$.
 
 \item
 Its Duggal transform ${T}^{D}=\widetilde {T}(1) $ is complex symmetric if and only if $ |\lambda_i|= |\lambda_{n-1-i}|$ for evey $1\leq i\leq n-2$.
 \end{enumerate}
\end{theorem}

{\it Proof.} 

As in the previous section, one may check easily
that

\[
\begin{array}{l}
\widetilde {T}(t)=|T|^{t}U |T|^{1-t}\\
 = \begin{pmatrix}
0 &   &  &  &  \\
 & |\lambda_{1}| &  &  &  \\
 &  &  &  \ddots &  \\
 &  &  &  &  |\lambda_{n-1}| \\
\end{pmatrix} ^t\begin{pmatrix}
0 & 1 & 0 & & \ldots&  0\\
0 & 0 & 1 & 0 & \ldots & 0\\
\vdots &\vdots &\ddots &\ddots & \ldots& 0\\
\vdots & \vdots & . & 0 & \ddots &  0\\
.& . & . & . & 0 & 1\\
0 & 0 & . & . &\ldots & 0 \\
\end{pmatrix}
\begin{pmatrix}
0 &   &  &  &  \\
 & |\lambda_{1}| &  &  &  \\
 &  &  &  \ddots &  \\
 &  &  &  &  |\lambda_{n-1}| \\
\end{pmatrix} ^{1-t}\\
= \begin{pmatrix}
0 & 0 &  &  &  &  \\
 & 0 & |\lambda_{1}|^{t}|\lambda_{2}|^{1-t} &  &  &  \\
  &  & 0 &  &  &  \\
 &  &  & 0 &  \ddots &  \\
 &  &  &  &   \ddots &|\lambda_{n-2}|^{t} |\lambda_{n-1}|^{1-t}\\
  &  &  &  &  &  0\\
\end{pmatrix}
=0\oplus
\begin{pmatrix}
0 &  |\lambda_{1}|^{t}|\lambda_{2}|^{1-t}     & 0 & \ldots & 0\\
\vdots &\ddots &\ddots & \ldots& 0\\
 \vdots & . & 0 & \ddots &  0\\
 . & . & . & 0 &   |\lambda_{n-2}|^{t} |\lambda_{n-1}|^{1-t}      \\
 0 & . & . &\ldots & 0 \\
\end{pmatrix}
\end{array}
\]
for every $t\in ]0,1]$, and the result follows immediately from \cite[Lemma 1]{GW1} and \cite[Theorem 3.1]{LZ2}.

\medskip

\begin{remark}\label{R2}

\begin{enumerate}
 \item 
If an operator $T$ is complex symmetric then its Aluthge transform  $\widetilde {T}=\widetilde {T}(\frac{1}{2}) $ is complex symmetric but the converse is not true: Consider for example $n=5$  and $|\lambda_{1}|=|\lambda_{3}|$,    $|\lambda_{2}|=|\lambda_{4}|$ and $|\lambda_{1}|\neq |\lambda_{4}|$.
 \item
 The second assertion of the theorem shows that for most  cases in this situation, the Duggal  transform is not complex symmetric.
 
 The explanation of the confusion in \cite{Ga} comes from the following: 
 as it is stated in 
\cite[Theorem 2]{Ga 4}, 
if $T\in{\mathcal{L}}(\mathcal{H})$ is a complex symmetric operator with a conjugation $C$
then there exists a partial conjugation $J$ supported on $\overline{ran (|T|)}$ such that $T=CJ|T|$ and $J|T|=|T|J.$
 A generalization of a theorem of Godi\v{c} and Lucenko is used to show  that the $U$ appearing in the polar decomposition may  be written as $U=CJ$ where $J$
 is partial conjugation which can of course be extended to a conjugation (let's say $\JJ$) acting on the whole space $\mathcal{H}$ 
 without affecting $T=CJ|T|=C{\JJ}|T|$. The only problem is that if one considers 
 $|T| C{\JJ}$, then it is not necessarily the Duggal transform of $T$.
 
 \medskip
 
 Indeed,
 
 Consider our previous example in Remark \ref{R1}

\[
T= \begin{pmatrix}
0 &  1 & 0 &  0 \\
0 &  0 & 2 &  0 \\
0 &  0 & 0 &  1 \\
0 &  0& 0 &  0 
 \end{pmatrix}=U|T|
\]
 
 where

\[
U= \begin{pmatrix}
0 &  1 & 0 &  0 \\
0 &  0 & 1 &  0 \\
0 &  0 & 0 &  1 \\
0 &  0& 0 &  0 
 \end{pmatrix}  \ \ \ {\text and } \ \ \ |T|= \begin{pmatrix}
0&  0 & 0 &  0 \\
0 &  1 & 0 &  0 \\
0 &  0 & 2 &  0 \\
0 &  0& 0 &  1 
 \end{pmatrix}.
\]

We know that $T$ is complex symmetric operator with the conjugation 
$C(z_1,z_2,z_3, z_4)= (\bar z_4, \bar z_3,  \bar z_2,\bar z_1)$.
We know also that $U=CJ$. Thus $J=CU$ and we have
$J(z_1,z_2,z_3, z_4)= ( 0, \bar z_4,  \bar z_3,\bar z_2)$.
(or equivalently $J$ is a partial conjugation such that $Je_1=0$, $Je_2=e_4 $, $Je_3=e_3 $ and $Je_4=e_2 $)

Obviously, $J$ can be extended to a conjugation ${\JJ}$ (by setting ${\JJ}e_1=e_1$
which means ${\JJ}(z_1,z_2,z_3, z_4)= ( \bar z_1,  \bar z_4,\bar z_3, \bar z_2)$.

It's rather easy to see that

\[
CJ= \begin{pmatrix}
0 &  1 & 0 &  0 \\
0 &  0 & 1 &  0 \\
0 &  0 & 0 &  1 \\
0 &  0& 0 &  0 
 \end{pmatrix}=U
 \ \ \ {\text while} \ \
C{\JJ}= \begin{pmatrix}
0 &  1 & 0 &  0\\
0 &  0 & 1 &  0 \\
0 &  0& 0 &  1 \\
1 &  0 & 0 &  0 
 \end{pmatrix}
\]

\[
|T|C{\JJ}= \begin{pmatrix}
0 &  0 & 0 &  0\\
0 &  0 & 1 &  0 \\
0 &  0& 0 &  2 \\
1 &  0 & 0 &  0 
 \end{pmatrix} \neq 
 \begin{pmatrix}
0 &  0 & 0 &  0\\
0 &  0& 1&  0 \\
0 &  0& 0 &  2 \\
0 &  0 & 0 &  0
 \end{pmatrix}=T^{D}.
\]

\item When $n=3$, all Duggal transforms of our studied operators are complex symmetric. Indeed, in this case, the Duggal transforms are nilpotent of degree 2 and it is known that these operators are complex symmetric.

 \end{enumerate}
\end{remark}

\section{Added remarks on mean transforms and the symmetric property}

 Recall that if  $T\in{\LL(\HH)}$, then the generalized mean transform of $T$ is the operator $\widehat{T}(t)=\frac{1}{2}[\widetilde {T}(t)+\widetilde{T}(1-t)]$,
where $\widetilde {T}(t)=|T|^{t}U |T|^{1-t}$ for $t\in (0,\frac{1}{2})$ is the generalized Aluthge transform of $T$.
The mean transform has been considered in \cite{LLY}.

\smallskip

%One may see the following result as a generalization of the one given in section xxxxx

\begin{theorem}\label{mean}
If $T= \sum_{i=1}^{n-1} \lambda_i e_{i}\otimes e_{i+1}$  and $\lambda_i\ne 0$ for all $i$, then
 its  generalized mean transform $\widehat {T}(t)$  ($t\in ]0,\frac{1}{2}]$)  is complex symmetric if and only if 
 $ |\lambda_i|^{t}|\lambda_{i+1}|^{1-t} +  |\lambda_{i}|^{1-t}|\lambda_{i+1}|^{t}= |\lambda_{n-1-i}|^{t}  |\lambda_{n-i}|^{1-t}+|\lambda_{n-1-i}|^{1-t}  |\lambda_{n-i}|^{t}
 $ for evey $1\leq i\leq n-2$ .
 In particular
 
 \begin{enumerate}
 \item 
 If $T$ is  complex symmetric then  its generalized mean transforms $  \widehat {T}(t)$ are complex symmetric 
 for all $t$ in $]0,\frac{1}{2}]$.

 \item
 On the other hand, its mean transform $\widehat {T}= \widehat {T}(0)$ is not complex symmetric in general.

 \end{enumerate}
\end{theorem}

{\it Proof.} 

As in the proof of theorem \ref{gAlu}, one has for $t\in ]0,\frac{1}{2}]$ that

\[
\begin{array}{l}
\widehat{T}(t)=\frac{1}{2}[\widetilde {T}(t)+\widetilde{T}(1-t)]\\

= \begin{pmatrix}
0 & 0 &  &  &  &  \\
 & 0 & \frac{1}{2}\Big[|\lambda_{1}|^{t}|\lambda_{2}|^{1-t}  +|\lambda_{1}|^{1-t}|\lambda_{2}|^{t} \Big] &  &  &  \\
  &  & 0 &  &  &  \\
 &  &  & 0 &  \ddots &  \\
 &  &  &  &   \ddots &\frac{1}{2} \Big[|\lambda_{n-2}|^{t} |\lambda_{n-1}|^{1-t}  + |\lambda_{n-2}|^{1-t} |\lambda_{n-1}|^{t} \Big] \\
  &  &  &  &  &  0\\
\end{pmatrix}\\
=0\oplus
\frac{1}{2}\begin{pmatrix}
0 &  \Big(|\lambda_{1}|^{t}|\lambda_{2}|^{1-t}  +|\lambda_{1}|^{1-t}|\lambda_{2}|^{t} \Big)   & 0 & \ldots & 0\\
\vdots &\ddots &\ddots & \ldots& 0\\
 \vdots & . & 0 & \ddots &  0\\
 . & . & . & 0 & \Big( |\lambda_{n-2}|^{t} |\lambda_{n-1}|^{1-t}  + |\lambda_{n-2}|^{1-t} |\lambda_{n-1}|^{t} \Big)       \\
 0 & . & . &\ldots & 0 \\
\end{pmatrix},
\end{array}
\]

which proves, thanks to \cite[Lemma 1]{GW1} and \cite[Theorem 3.1]{LZ2}, the main statement of our theorem.
The second statement is obvious since  $|\lambda_{i}|=|\lambda_{n-i}|$ and $|\lambda_{i+1}|=|\lambda_{n-i-1}|$ for every $1\leq i\leq n-2$, whenever $T$
is complex symmetric.

\smallskip

Now, the last statement is illustrated by the following example.

Let's consider the example with $n=4$. (One could, of course, treat the general case; we choose to  leave it to the interested reader.
)

\[
T= \begin{pmatrix}
0 &   |\lambda_{1}|& 0 &  0 \\
0 &  0 & |\lambda_{2}| &  0 \\
0 &  0 & 0 & |\lambda_{3}| \\
0 &  0& 0 &  0 
 \end{pmatrix}
\]

It's rather easy to show that 

\[
\widehat {T}=
 \begin{pmatrix}
0 &   \frac{ |\lambda_{1}|}{ 2}& 0 &  0 \\
0 &  0 &   \frac{ |\lambda_{1}|+|\lambda_{2}| }{2} &  0 \\
0 &  0 & 0 &  \frac{ |\lambda_{2}|+|\lambda_{3}| }{2}    \\
0 &  0& 0 &  0 
 \end{pmatrix}
\]

$\bullet$  Then $\widehat {T}$ is complex symmetric if and only if $|\lambda_{1}|=|\lambda_{2}|+|\lambda_{3}|$.
Notice that this happens only if $\lambda_{2}=0$ if we assume, in addition, that $T$ is complex symmetric (which means $|\lambda_{1}|=|\lambda_{3}|$).

$\bullet$ $\widehat {T}$ may be complex symmetric even though $T$ is not!

$\bullet$ Better (or worse!), for our examples  $\widehat {T}$ is never complex symmetric when $T$ and its Duggal transform $T^D$ are!

Indeed, it is enough to check it for the following example.

If

\[
T= \begin{pmatrix}
0 &  1 & 0 &  0 \\
0 &  0 & 1 &  0 \\
0 &  0 & 0 &  1 \\
0 &  0& 0 &  0 
 \end{pmatrix}
\ \ \ \textit{then} \ \ \ \
\widehat {T}=
\begin{pmatrix}
0 &  \frac{1}{2} & 0 &  0 \\
0 &  0 & 1 &  0 \\
0 &  0 & 0 &  1 \\
0 &  0& 0 &  0 
 \end{pmatrix}
\]

and $\widehat {T}$ is not  complex symmetric.

\begin{remark}\label{R3}
 It's easy to see that the considered examples:
$T= \sum_{i=1}^{n-1} \lambda_i e_{i}\otimes e_{i+1}$  and $\lambda_i\ne 0$ for all $i$, 
are also centered in the sense that the doubly infinite sequence (here it is a finite sequence, since T is nilpotent $T^n=0$)

$$ \{ \dots, (T^2)^* T^2, T^*T, TT^*, T^2(T^2)^*, \dots \}$$
is a set of mutually commuting operators.

This  also answers in the negative some questions asked in section 5 of \cite{HMT}.
\end{remark}

\vspace{5mm}

\end{document}